\documentclass[12pt]{article}
\usepackage{amsfonts,amsmath,amsthm,yfonts,mathdots,
latexsym,graphicx,amssymb,hyperref,url,fancyhdr}

\allowdisplaybreaks
\newtheorem{theorem}{Theorem}[section]
\newtheorem{lemma}[theorem]{Lemma}

\newtheorem{conjecture}{Conjecture}

\newtheorem{question}{Question}
\title{Triangle-free graphs with the maximum number of cycles}
\author{Andrii Arman\thanks{\url{armana@cc.umanitoba.ca}}
\and  David S. Gunderson\thanks{ \url{David.Gunderson@umanitoba.ca}.
  Partially supported by NSERC grant 228064-00.}
\and Sergei Tsaturian\thanks{\url{tsaturis@cc.umanitoba.ca}} \\
  \footnotesize University of Manitoba, Winnipeg, Manitoba, Canada, R3T 2N2. }
\date{22 January 2015}
\begin{document}
\maketitle
\begin{abstract}
It is shown that for $n\geq 141$, among all triangle-free graphs on  $n$ vertices,
the  balanced complete bipartite graph $K_{\lceil n/2 \rceil, \lfloor n/2 \rfloor}$
is the unique triangle-free graph with the
maximum number of cycles.  Using modified Bessel functions, tight estimates are
 given for the  number of cycles in  $K_{\lceil n/2 \rceil, \lfloor n/2 \rfloor}$.
 Also, an upper bound for the number of Hamiltonian cycles in a triangle-free
 graph is given.
\end{abstract}

\section{Introduction}

In a  recent article \cite{DGLS:15},  the maximum number of cycles
in a triangle-free graph was considered.  It was asked which triangle-free graphs contain the maximum  number of cycles; this question arose from the study of path-finding algorithms \cite{BCD:13}.
The same authors posed the following conjecture:
\begin{conjecture}[Durocher--Gunderson--Li--Skala, 2014 \cite{DGLS:15}]\label{con:main}
For each $n\geq 4$,  the balanced complete bipartite graph
$K_{\lceil n/2\rceil, \lfloor n/2 \rfloor}$ contains more cycles than any other
$n$-vertex triangle-free graph.
\end{conjecture}
The authors \cite{DGLS:15} confirmed Conjecture \ref{con:main} when $4\leq n\leq 13$, and made progress toward this conjecture in general.  For example, they showed the conjecture
to be true when restricted to ``nearly regular graphs'', that is, for
each positive integer $k$ and sufficiently large $n$, $K_{\lceil n/2\rceil, \lfloor n/2 \rfloor}$ has
more cycles than any other triangle-free graph on $n$ vertices whose minimum degree
and maximum degree differ by at most $k$.

 In Theorems \ref{th:main}  and \ref{th:mainngeq141} below,
 it is shown that Conjecture \ref{con:main} holds true for $n\geq 141$.
 Theorem \ref{th:cTn2approx} gives a useful estimate for
 the number of cycles in $K_{\lceil n/2\rceil, \lfloor n/2 \rfloor}$. In Lemma
 \ref{le:2ham}, an upper bound is given for the number of Hamiltonian cycles in
 a triangle-free graph.

Even though  Conjecture \ref{con:main} arose from a very specific problem in
computing, it can be considered as a significant problem in two aspects of graph
theory: counting cycles in graphs, and the structure of
triangle-free graphs.
In recent decades, bounds have been proved for the maximum number of cycles in
various classes of graphs.  Some of these classes include
complete graphs ({\it e.g.} \cite{HaMa:71}),
planar graphs  ({\it e.g.}, \cite{AT:08}, \cite{AFK:99}, \cite{BKKSS:07}),
outerplaner graphs  and series-parallel graphs ({\it e.g.}, \cite{dMN:12}),
graphs with large maximum degree ({\it e.g.}, \cite{BBRS:03} for those without a
specified odd cycle),
graphs with specified minimum degree (see, {\it e.g.}, \cite{Volk:96})
graphs with a specified cyclomatic number or number of  edges({\it e.g.}, \cite{Ahrens:97},
 \cite{Kirc:47}, \cite{EnSl:81}, \cite{Guic:96}; see also \cite[Ch4, Ch10]{Konig:35}),
cubic graphs ({\it e.g.}, \cite{AT:97}), \cite{CT:12}),
graphs with fixed girth ({\it e.g.} \cite{Mark:04}),
$k$-connected graphs ({\it e.g.} \cite{Knor:94}),
Hamiltonian graphs ({\it e.g.}, \cite{RaSt:05}, \cite{Shi:94}, \cite{Volk:96})
Hamiltonian graphs with a fixed number of edges ({\it e.g.} \cite{GuSh:93}),
2-factors of the de Bruijn graph ({\it e.g.}, \cite{Fred:70}),
graphs with a cut-vertex ({\it e.g.}, \cite{Volk:96})
complements of trees ({\it e.g.}, \cite{HoWi:85}, \cite{Reid:76},  \cite{Zhou:88}),
and random graphs ({\it e.g.}, \cite{Taka:88}).
In some cases, the structure of the extremal graphs are also found
({\it e.g.}, \cite{BBRS:03}, \cite{RaSt:05}).

In 1973, Erd\H{o}s, Kleitman, and Rothschild \cite{EKR:76} showed that
for $r\geq 3$, as $n\rightarrow\infty$,
the number of $K_r$-free graphs on $n$ vertices is
\[2^{(1-\frac{1}{r-1} +o(1))\binom{n}{2}}.\]
As a consequence, the number
of triangle-free graphs is very close to the number of bipartite graphs,
and so almost all triangle-free graphs are bipartite.
By Mantel's theorem \cite{Mantel:07}, among graphs on $n$ vertices,
the triangle-free graph with
the most number of edges is the balanced complete bipartite graph
$K_{\lceil n/2 \rceil, \lfloor n/2 \rfloor}$.

Since $K_{\lceil n/2 \rceil, \lfloor n/2 \rfloor}$ is the triangle-free graph on $n$ vertices with the
most number of edges, and nearly all triangle-free graphs are bipartite,
Conjecture \ref{con:main}  might seem reasonable, even though
$K_{\lceil n/2 \rceil, \lfloor n/2 \rfloor}$ contains no odd cycles.

\section{Notation and approximations used}\label{se:notappns}
A graph $G$ is an ordered pair $G=(V,E)=(V(G),E(G))$,
where $V$ is a non-empty set and $E$ is a set of unordered pairs from $V$. Elements of $V$ are called vertices and elements of $E$ are called edges. Under this definition,
graphs are simple, that is, there are no loops nor multiple edges.
% The complete graph
%on $n$ vertices is denoted by $K_n$ (where every one of the possible $\binom{n}{2}$ pairs
%of vertices is an edge).

 An edge $\{x,y\}\in E(G)$  is denoted by simply $xy$.
 The neighbourhood of any vertex $v\in V(G)$ is
 $N_G(x)=\{y\in V(G):xy\in E(G)\}$, and the degree of $x$ is $\deg_G(x)=|N(x)|$.
When it is clear what $G$ is, subscripts are deleted, using only $N(x)$ and $\deg(x)$.
The minimum degree of vertices in a graph $G$ is denoted by $\delta(G)$, and the maximum degree is denoted $\Delta(G)$.
If $Y\subset V(G)$, the subgraph of $G$ induced by $Y$ is denoted $G[Y]$.

A graph $G=(V,E)$ is called bipartite iff there is a partition $V=A\cup B$ so that
$E\subset \{\{x,y\}:x\in A, y\in B\}$; if $E=\{\{x,y\}:x\in A, y\in B\}$, then $G$ is called
the complete bipartite graph on partite sets $A$ and $B$, denoted $G=K_{|A|,|B|}$.
The balanced complete bipartite graph on $n$ vertices is
$K_{\lfloor n/2\rfloor, \lceil n/2\rceil}$.
%this graph is also called a Tur\'an graph, denoted by $T(n,2)$
A cycle on $m$ vertices is denoted $C_m$. The complement
of a graph $G$ is denoted $\overline{G}$.
For any graph $G$, let $c(G)$ denote the number of cycles in $G$.

The number $e$ is the base of the natural logarithm.
Stirling's approximation formula says that as $n\rightarrow \infty$,
\begin{equation}\label{eq:stirlingsformula}
n!=(1+o(1))\sqrt{2\pi n}(n/e)^n.
\end{equation}
In 1955, Robbins \cite{Robbins:55} proved that
\[\sqrt{2\pi n}\left(\frac{n}{e} \right)^ne^{\frac{1}{12n+1}}
<n!<\sqrt{2\pi n}\left(\frac{n}{e} \right)^ne^{\frac{1}{12n}}.\]
Slightly more convenient bounds are used (that are valid for all $n\geq 1$):
\begin{equation}\label{eq:simplebdsn!}
%n!< \frac{n^{n+\frac{1}{2}}}{e^{n-1}}.
\sqrt{2\pi}\cdot \sqrt{n}\left(\frac{n}{e}\right)^n < n! < e\cdot \sqrt{n}\left(\frac{n}{e}\right)^n.
\end{equation}
%Another simple inequality says that for $x>0$,
%\begin{equation}\label{eq:boundetox}
% 1+x\leq e^x.
%\end{equation}
%
%\begin{lemma} Let $e$ be the base for the natural logarithm. For $k=o(n^{3/4})$,
%\begin{equation}\label{eq:nchoosekapprox}
%\binom{n}{k}\sim \frac{1}{\sqrt{2\pi k}}\left(\frac{ne}{k}\right)^k.
%\end{equation}
%\end{lemma}
%
%In many calculations, the following simple bounds are sufficient:
%\begin{lemma}\label{le:binombounds}
%For any positive integers $k\leq n$,
%\[\left(\frac{n}{k}\right)^k \leq \binom{n}{k} \leq \left(\frac{en}{k}\right)^k.\]
%\end{lemma}

Two modified Bessel functions (see, {\it e.g.}, \cite{AbSt:72}) are used:
\begin{equation}\label{eq:I_0x}
I_0(x)= \sum_{k=0}^{\infty}\frac{x^{2k}}{2^{2k}(i!)^2}.
\end{equation}
\begin{equation}\label{eq:I_1x}
I_1(x)= \sum_{k=0}^{\infty}\frac{x^{2k+1}}{2^{2k+1}i!(i+1)!}.
\end{equation}
 In particular, when $x=2$ is used in either modified Bessel function,
 useful approximations are obtained:

 \begin{equation}\label{eq:I_02}
2.27958 \leq \sum_{k=0}^{\infty}\frac{1}{(i!)^2}=I_0(2)  \leq 2.279586;
 \end{equation}

 \begin{equation}\label{eq:I_12}
1.5906 \leq \sum_{i=0}^{\infty}\frac{i}{(i!)^2}=\sum_{k=0}^{\infty}\frac{1}{k!(k+1)!}=I_1(2) \leq 1.59064.
 \end{equation}
% and, with combining (\ref{eq:I_02}) and (\ref{eq:I_12}),
% \begin{equation}\label{eq:fact3}
% \sum_{i=0}^{\infty}\frac{(2i-1)i}{(i!)^2}=\sum_{i=1}^{\infty}\frac{2}{((i-1)!)^2}-\sum_{i=0}^{\infty}\frac{i}{(i!)^2}=2I_0(2)-I_1(2)\sim 2.968.
% \end{equation}
%

\section{Preliminaries}
%
%Say that a graph $G$ on $n$ vertices is
%{\em cycle-maximal} triangle-free iff $G$ has the greatest number of cycles among
%all triangle-free graphs on $n$ vertices.
%\begin{lemma}[\cite{DGLS:15}]
% If $G$ is a cycle-maximal triangle-free graph on at least 4 vertices,
%then $G$ is 2-connected.
%\end{lemma}
%
%\begin{lemma}[\cite{DGLS:15}]
%Let $G$ be a cycle-maximal triangle-free graph with at least 4 vertices.
%Then every edge of $G$ is some 4-cycle, and the addition of any missing
%edge forms a triangle.
%\end{lemma}
%If one restricts to ``nearly regular'' graphs, Conjecture \ref{con:main} is true for $n$ sufficiently
%large:
%
%\begin{lemma}[\cite{DGLS:15}]
%For each positive integer $k$, there exists $n_k$ so that for $n\geq n_k$,
%among all graphs with $\Delta(G)-\delta(G)\leq k$,  the graph
%$K_{\lceil n/2\rceil, \lfloor n/2 \rfloor}$ contains the
%most number of cycles.
%\end{lemma}

The following shows that among all bipartite graphs, the balanced one has
the most cycles.
\begin{lemma}[\cite{DGLS:15}]\label{le:maxcbipartite}
For $n\geq 4$, among all bipartite graphs on $n$ vertices,
$K_{\lceil n/2\rceil, \lfloor n/2 \rfloor}$ has the greatest number of cycles; that is,
$K_{\lceil n/2\rceil, \lfloor n/2 \rfloor}$ is the unique cycle-maximal bipartite graph on $n$ vertices.
\end{lemma}

So, to settle Conjecture \ref{con:main}, it is then sufficient to prove that a cycle-maximal triangle-free
graph is bipartite. To this end, the following result is essential:
\begin{theorem}[Andr\'asfai, 1964 \cite{Andr:64}]\label{th:andr64}
 Any triangle-free graph $G$  on $n$ vertices with $\delta(G)> 2n/5$
 is bipartite.
 \end{theorem}
See also \cite{AES:74} for an English proof of Theorem \ref{th:andr64} and related results.
 Theorem \ref{th:andr64} is sharp because of $C_5$ (or a blow-up of
 $C_5$).

\begin{lemma}[\cite{DGLS:15}]\label{le:cTn2}
For $n\geq 4$, the number of cycles in the balanced complete bipartite graph is
\begin{equation}\label{eq:cTn2}
c(K_{\lfloor n/2\rfloor, \lceil n/2\rceil})=\sum_{k=2}^{\lfloor n/2\rfloor}
\frac{\lfloor n/2\rfloor ! \lceil n/2\rceil !}{2k(\lfloor n/2 \rfloor -k)! (\lceil n/2 \rceil -k)!}.
\end{equation}
\end{lemma}

The following form for the number of cycles in $K_{\lfloor n/2\rfloor, \lceil n/2\rceil}$
gives a way to estimate the right hand side of (\ref{eq:cTn2}) in Lemma \ref{le:cTn2}:
\begin{theorem}\label{th:cTn2approx}
For $n\geq 12$,
 \begin{align}\label{eq:cTn2lowerbd}
 c(K_{\lfloor n/2\rfloor, \lceil n/2\rceil})&\geq
 \frac{\lfloor n/2\rfloor ! \lceil n/2\rceil !}{2\lfloor n/2 \rfloor} \cdot  \begin{cases}
 I_0(2) \cdot  &\mbox{if $n$ is even}\\
 I_1(2)  \cdot  &\mbox{if $n$ is odd}.
\end{cases}\\
&\geq \pi \left(\frac{n}{2e}\right)^n\cdot \begin{cases}
 I_0(2)  &\mbox{if $n$ is even}\\
 I_1(2)  &\mbox{if $n$ is odd}.
\end{cases}
\label{eq:cTn2lowerbdsimple}
\end{align}
and as $n\rightarrow \infty$,
\begin{equation}\label{eq:cTn2approx}
c(K_{\lfloor n/2\rfloor, \lceil n/2\rceil})=(1+o(1))
\begin{cases}
 I_0(2)\pi \left(\frac{n}{2e}\right)^n &\mbox{if $n$ is even}\\
 I_1(2)  \pi \left(\frac{n}{2e}\right)^n&\mbox{if $n$ is odd}
 \end{cases}.
 \end{equation}
 \end{theorem}

\noindent{\bf Proof:} Using (\ref{eq:simplebdsn!}), the proof that (\ref{eq:cTn2lowerbdsimple})
follows from (\ref{eq:cTn2lowerbd}) is elementary,  and so is omitted.

By Lemma \ref{le:cTn2}, write
\begin{align}
c(K_{\lfloor n/2\rfloor, \lceil n/2\rceil})
&=\sum_{k=2}^{\lfloor n/2\rfloor}
\frac{\lfloor n/2\rfloor ! \lceil n/2\rceil !}{2k(\lfloor n/2 \rfloor -k)! (\lceil n/2 \rceil -k)!} \notag\\
&=\frac{\lfloor n/2\rfloor ! \lceil n/2\rceil !}{2\lfloor n/2 \rfloor}      \cdot
   \sum_{k=2}^{\lfloor n/2\rfloor}\frac{\lfloor n/2\rfloor }{k(\lfloor n/2 \rfloor -k)! (\lceil n/2 \rceil -k)!}.
   \label{eq:cTn2twofracs}
\end{align}

Case 1 ($n$ even): Suppose that for $\ell \geq 2$,  $n=2\ell$, and set
\[a_\ell=\sum_{k=2}^\ell \frac{\ell}{k((\ell-k)!)^2}=\sum_{i=0}^{\ell-2}\frac{\ell}{(\ell-i)(i!)^2}.\]
Then equation (\ref{eq:cTn2twofracs}) becomes
\begin{equation}\label{eq:cT2nfracdotan}
c(K_{\lfloor n/2\rfloor, \lceil n/2\rceil})
%=\sum_{k=2}^{\lfloor n/2\rfloor}
%\frac{\lfloor n/2\rfloor ! \lceil n/2\rceil !}{2k(\lfloor n/2 \rfloor -k)! (\lceil n/2 \rceil -k)!}
=\frac{\lfloor n/2\rfloor ! \lceil n/2\rceil !}{2\lfloor n/2 \rfloor}\cdot a_\ell.
\end{equation}
\medskip

\noindent Claim: For $\ell \geq 6$, $a_{\ell+1} \leq a_{\ell}$.

\noindent Proof of Claim:
\begin{align*}
a_{\ell} - a_{\ell+1}
& = \sum_{i=0}^{\ell-2} \left(\frac{\ell}{\ell-i}-\frac{\ell+1}{\ell+1-i}\right) \frac{1}{(i!)^2}
             -\frac{\ell+1}{2((\ell-1)!)^2}\\
&= \sum_{i=0}^{\ell-2} \left(\frac{i}{(\ell+1)(\ell-i)}\right)\frac{1}{(i!)^2}
             -\frac{\ell+1}{2((\ell-1)!)^2}\\
&= \sum_{i=2}^{\ell-2} \frac{i}{(\ell+1-i)(\ell-i)}\cdot \frac{1}{(i!)^2} -\frac{1}{\ell(\ell-1)}
                   -\frac{\ell+1}{2((\ell-1)!)^2}\\
&\geq 0+\frac{2((\ell-1)!)^2-(\ell+1)\ell(\ell-1)}{2((\ell-1)!)^2}\\
&\geq 0 &&\mbox{(for $\ell\geq 6$)},
\end{align*}
finishing the proof of the claim.
\medskip

Since the sequence $\{a_{\ell}\}$ is non-increasing and bounded below (by 0, {\it e.g.}),
$\lim_{\ell\rightarrow\infty}a_{\ell}$ exists. To find this limit, first apply partial fractions:
\[a_\ell=\sum_{i=0}^{\ell-2}\frac{\ell}{(\ell-i)(i!)^2}
=\sum_{i=0}^{\ell-2}\frac{1}{(i!)^2}+\sum_{i=0}^{\ell-2}\frac{i}{(\ell-i)(i!)^2}.\]
Put $\displaystyle b_\ell=\sum_{i=0}^{\ell-2}\frac{1}{(i!)^2}$ and
      $\displaystyle c_\ell = \sum_{i=0}^{\ell-2}\frac{i}{(\ell-i)(i!)^2}$.
Then
\begin{align*}
c_\ell&=\sum_{i=0}^{\ell-2}\frac{i}{(\ell-i)(i!)^2}\\
&=\sum_{i=0}^{3}\frac{i}{(\ell-i)(i!)^2}+\sum_{i=4}^{\ell-2}\frac{i}{(\ell-i)(i!)^2}\\
&\leq \frac{3}{\ell-3}+\frac{1}{\ell}\sum_{i=4}^{\ell-2}\frac{1}{i!}&&\mbox{(since $\frac{i}{(\ell-i)i!}\leq \frac{1}{\ell}$ for $\ell\geq 4$})\\
&\leq \frac{3}{\ell-3}+\frac{e}{\ell},
\end{align*}
and therefore, $\lim_{\ell\rightarrow \infty}c_\ell=0$. Thus,
\begin{align*}
\lim_{\ell\rightarrow\infty}a_\ell &=\lim_{\ell\rightarrow\infty}(b_\ell+c_\ell)\\
&=\lim_{\ell\rightarrow\infty}b_\ell\\
&=\sum_{i=0}^{\infty}\frac{1}{(i!)^2}\\
&=I_0(2)&&\mbox{(by (\ref{eq:I_02}))}.
\end{align*}

Since $a_\ell$ is non-increasing for $\ell\geq 6$,  for $n\geq 12$,
\[c(K_{\lfloor n/2\rfloor, \lceil n/2\rceil})\geq \frac{\lfloor n/2\rfloor ! \lceil n/2\rceil !}{2\lfloor n/2 \rfloor}\cdot I_0(2),\]
which proves the even case of (\ref{eq:cTn2lowerbd}).
By (\ref{eq:I_02}), as $n\rightarrow \infty$,
\[c(K_{\lfloor n/2\rfloor, \lceil n/2\rceil}) =(1+o(1))\frac{\lfloor n/2\rfloor ! \lceil n/2\rceil !}{2\lfloor n/2 \rfloor}\cdot I_0(2),\]
and by Stirling's approximation (\ref{eq:stirlingsformula}), the proof of the even case of
 (\ref{eq:cTn2approx}) is complete.\medskip

Case 2 ($n$ odd): Suppose that for $\ell \geq 6$, $n=2\ell+1$. The proof
follows the even case, and so is only outlined.
Put
\[a_\ell=\sum_{k=2}^\ell\frac{\ell}{k(\ell-k)!(\ell+1-k)!}
=\sum_{i=0}^{\ell-2}\frac{\ell}{(\ell-i)i!(i+1)!}.\]

\noindent Claim: For $\ell\geq 4$, $a_{\ell+1}\leq a_\ell$.

\noindent{Proof of claim:} Letting $\ell \geq 4$,
\begin{align*}
a_\ell-a_{\ell+1}&=\sum_{i=0}^{\ell-2}\frac{i}{(\ell+1-i)(\ell-i)}\cdot \frac{1}{i!(i+1)!}
                                                                           -\frac{\ell+1}{2(\ell-1)!\ell !}\\
&=\sum_{i=2}^{\ell-2}\frac{i}{(\ell+1-i)(\ell-i)}\cdot \frac{1}{i!(i+1)!} - \frac{1}{2(\ell-1)\ell}
                                                                           -\frac{\ell+1}{2(\ell-1)!\ell !}\\
&\geq 0+\frac{(\ell-2)!(\ell-1)!-(\ell+1)}{2(\ell-1)!\ell !}\\
&\geq 0,
\end{align*}
finishing the proof of the claim.

Therefore, $\lim_{\ell\rightarrow\infty}a_\ell$ exists. To find this limit, write
\[a_\ell=\sum_{i=0}^{\ell-2}\frac{1}{i!(i+1)!}+ \sum_{i=0}^{\ell-2}\frac{i}{(\ell-i)i!(i+1)!}.\]
Letting $\displaystyle b_\ell =\sum_{i=0}^{\ell-2}\frac{1}{i!(i+1)!}$ and
  $\displaystyle c_{\ell}= \sum_{i=0}^{\ell-2}\frac{i}{(\ell-i)i!(i+1)!}$, observe that
  \[  c_\ell= \sum_{i=0}^{\ell-2}\frac{i}{(\ell-i)i!(i+1)!}+ \sum_{i=0}^{\ell-2}\frac{i}{(\ell-i)i!(i+1)!}
           \leq \frac{3}{\ell-3}+\frac{e}{\ell},\]
 and so $\displaystyle \lim_{\ell\rightarrow\infty}c_\ell=0$. Thus,
 \[\lim_{\ell\rightarrow\infty}a_\ell=\lim_{\ell\rightarrow\infty}b_\ell
    =\sum_{i=0}^{\infty}\frac{1}{i!(i+1)!}
    =\sum_{i=0}^\infty \frac{i+1}{((i+1)!)^2}
    =\sum_{i=0}^\infty \frac{i}{(i!)^2},\]
    which, by (\ref{eq:I_12}), is equal to $I_1(2)$.
    Then again
  \begin{align*}
  c(K_{\lfloor n/2\rfloor, \lceil n/2\rceil})
   & \geq \frac{\lfloor n/2\rfloor ! \lceil n/2\rceil !}{2\lfloor n/2 \rfloor}  \cdot I_1(2)\\
   & = \frac{ \ell ! (\ell+1) !}{2\ell}  \cdot I_1(2)\\
   &=       \frac{ (\ell !)^2}{2\ell}(\ell+1) \cdot I_1(2)\\
   & =(1+o(1))\pi\left(\frac{\ell}{e}\right)^{2\ell}(\ell+1)\cdot I_1(2)&&\mbox{(by (\ref{eq:stirlingsformula}))}\\
    & >(1+o(1))\pi\left(\frac{\ell}{e}\right)^{2\ell}(\ell-1)\cdot I_1(2)\\
   & =1+o(1))\pi \left(\frac{n-1}{2e}\right)^{n-1} \left(\frac{n-1}{2}\right)  \cdot I_1(2)\\
   &=(1+o(1))\pi e\left(\frac{n-1}{2e}\right)^{n}\cdot I_1(2)\\
   &=(1+o(1))\pi e\left(\frac{n-1}{n}\right)^n\left(\frac{n}{2e}\right)^{n}\cdot I_1(2)\\
    &=(1+o(1))\pi \left(\frac{n}{2e}\right)^n  \cdot I_1(2),
 \end{align*}
 and as $n\rightarrow\infty$,
  \[c(K_{\lfloor n/2\rfloor, \lceil n/2\rceil})=(1+o(1))\pi \left(\frac{n}{2e}\right)^n  \cdot I_1(2).\]
  This completes the proof for odd $n$, and so the proof of the lemma.
  \qed\bigskip

\begin{lemma}\label{le:6vcxy}
 Let $H$ be a triangle-free graph on 6 vertices with $x,y\in V(H)$. Then
there are at most 9 different $x$--$y$ paths.
\end{lemma}
\noindent{\bf Proof:} Consider two cases.

Case 1: $H$ contains no copy of $C_5$. Then $H$ contains no odd cycle,
and so is bipartite. Without loss of generality, add edges to $H$ to make $H$
a complete bipartite graph. There are only four different complete bipartite
graphs on six vertices, namely $\overline{K_6}$, $K_{1,5}$, $K_{2,4}$, and $K_{3,3}$.
By inspection, in any of these, the maximum number of paths between any two
 vertices is at most 9.

Case 2: $H$ contains a copy of $C_5$.  Suppose that $x_1,x_2,x_3,x_4,x_5,x_1$
forms a cycle $C$, and that $x_6$ is the remaining vertex. Then $x_6$ is adjacent to
at most two vertices of $C$. If $x_6$ is adjacent to fewer than two vertices of $C$, add
an extra edge or two so that $x_6$ is adjacent to precisely two vertices of $C$;
without loss, suppose that $x_6$ is adjacent to $x_1$ and $x_3$. Then the maximum
number of paths between any two vertices is 4 (for example, between $x_2$ and $x_6$).
\qed\bigskip

\section{Counting types of cycles}\label{se:countingtypes}
\begin{lemma}\label{le:1even}
There exists $n_0\in\mathbb{Z}^{+}$ so that for every even integer
$n\geq n_0$, if $G$ is a triangle-free graph on $n$ vertices, and
$x_1x_2\in E(G)$, then the number of cycles containing the edge $x_1x_2$
is at most $10\pi \frac{n^{n-1}}{(2e)^n}$.
\end{lemma}
\noindent{\bf Proof:} Let $G$ be a triangle-free graph on $n$ vertices,
and let $x_1x_2\in E(G)$.
For each $k=4,\ldots,n$, let $c_k$ denote the number of cycles
of length $k$ that contain the edge $x_1x_2$. The goal is to give an upper bound
for $\sum_{k=4}^nc_k$.

Let $2\leq i\leq \frac{n-4}{2}$;
 an upper bound on $c_{2i}+c_{2i+1}$ is first calculated; to do so, count all possible
 cycles of the form $x_1,x_2,\ldots,x_{2i}$ or $x_1,x_2,\ldots,x_{2i+1}$.
  For each $j>1$, there are at most $d_j=|N(x_j)\backslash\{x_1,\ldots,x_{j-1}\}|$ ways
   to choose an $x_{j+1}$. Note that $N(x_j)\cap N(x_{j+1})=\emptyset$,
 since otherwise a triangle is formed with $x_j$ and $x_{j+1}$. Also,
 \[|(N(x_j)\backslash\{x_1,\ldots,x_{j-1}\})\cup( N(x_{j+1})\backslash\{x_1,\ldots,x_{j}\})|
        \leq |V(G)\backslash\{x_1,\ldots,x_j\}|=n-j.\]
Therefore,
\begin{align*}
d_j+d_{j+1}
&\leq |N(x_j)\backslash\{x_1,\ldots,x_{j-1}\}|+| N(x_{j+1})\backslash\{x_1,\ldots,x_{j}\}|\\
&=[((N(x_j)\backslash\{x_1,\ldots,x_{j-1}\})\cup( N(x_{j+1})\backslash\{x_1,\ldots,x_{j}\})|\\
   &     \leq n-j,
  \end{align*}
  and thus
  \begin{equation}\label{eq:3star}
  d_jd_{j+1}\leq \left\lfloor \frac{n-j}{2}\right\rfloor \cdot \left\lceil \frac{n-j}{2}\right\rceil.
 \end{equation}
 Using  (\ref{eq:3star}), the number of ways to choose vertices
 $x_3,x_4,\ldots,x_{2i}$ so that $x_1,x_2,x_3,x_4,\ldots, x_{2i}$ form a path
 is at most
 \begin{equation}\label{eq:1star}
 \prod_{j=2}^{2i-1}d_j=\prod_{j=1}^{i-1}(d_{2j}d_{2j+1})
 \leq \prod_{j=1}^{i-1}\left(\left\lfloor \frac{n-2j}{2}\right\rfloor \cdot \left\lceil \frac{n-2j}{2}\right\rceil\right)
 =\prod_{j=1}^{i-1}\left(\frac{n-2j}{2}\right)^2.
 \end{equation}
 If there is an edge $x_{2i}x_1\in E(G)$, there is one cycle $x_1,x_2,\ldots,x_{2i}$
 of length $2i$, and no cycles of the form $x_1,x_2,\ldots,x_{2i+1}$ because
 otherwise, $x_1,x_{2i}, x_{2i+1}$ form a triangle. So, in total, there is exactly one
 cycle that contains the path $x_1,x_2,\ldots,x_{2i}$ and has length $2i$ or $2i+1$.
  If there is no edge $x_{2i}x_1$, there is no cycle $x_1,\ldots,x_{2i}$ and at most
 $n-2i$ cycles of the form $x_1,\ldots,x_{2i}x_{2i+1}$. In any case,
 there are at most $n-2i$ cycles containing the path $x_1,\ldots,x_{2i}$.

 By these observations and inequality (\ref{eq:1star}),
 \begin{equation}\label{eq:2star}
 c_{2i}+c_{2i+1}\leq (n-2i)\prod_{j=1}^{i-1}\left(\frac{n-2j}{2}\right)^2.
 \end{equation}
 To evaluate $\sum_{k=4}^nc_k$, separate the sum into two parts:
 \begin{align}
 \sum_{k=4}^{n-5}c_k
 &=\sum_{i=2}^{(n-6)/2}(c_{2i}+c_{2i+1})\notag \\
 &\leq \sum_{i=2}^{(n-6)/2}\left((n-2i)\prod_{j=1}^{i-1}\left(\frac{n-2j}{2}\right)^2\right)
                   &&\mbox{(by (\ref{eq:2star}))}\notag \\
 &=\sum_{i=2}^{(n-6)/2}(n-2i)\left(\frac{\left(\frac{n-2}{2}\right)!}{\left(\frac{n-2i}{2}\right)!} \right)^2
 \notag \\
 &=\left( \left(\frac{n-2}{2}\right)! \right)^2\sum_{j=3}^{\frac{n-4}{2}}\frac{2j}{(j!)^2}\notag\\
 &=\left( \left(\frac{n-2}{2}\right)! \right)^2\left(\sum_{j=1}^{\frac{n-4}{2}}\frac{2j}{(j!)^2} -\frac{2}{(1!)^2}-\frac{2\cdot 2}{(2!)^2}\right)\notag\\
 &\leq \left( \left(\frac{n-2}{2}\right)! \right)^2(2\cdot (1.591)-3)    &&\mbox{(by (\ref{eq:I_12}))} \notag\\
 &<0.19\left( \left(\frac{n-2}{2}\right)! \right)^2.\label{eq:5star}
 \end{align}
 To count $\displaystyle\sum_{k=n-4}^nc_k$, note that by (\ref{eq:3star}), there are at most
 \[\prod_{i=2}^{n-5}d_i\leq \prod_{j=1}^{\frac{n-6}{2}}\left(\frac{n-2j}{2}\right)^2\]
 ways to choose a path $x_1,x_2,\ldots,x_{n-4}$, and by Lemma \ref{le:6vcxy},
 there are at most 9 paths that connect $x_{n-4}$ and $x_1$ in the graph
 $G\backslash\{x_1,\ldots,x_{n-5}\}$; that is, there are at most 9 ways to complete
 the path $x_1,x_2,\ldots,x_{n-4}$ to a cycle. Therefore,
 \begin{equation}\label{eq:4star}
 \sum_{k=n-4}^nc_k\leq 9\prod_{j=1}^{\frac{n-6}{2}}\left(\frac{n-2j}{2}\right)^2
 =9\cdot \frac{\left(\left(\frac{n-2}{2}\right)! \right)^2}{(2!)^2}
 =\frac{9}{4}\left(\left(\frac{n-2}{2}\right)! \right)^2.
 \end{equation}
 Adding equations (\ref{eq:5star}) and (\ref{eq:4star}),
 \begin{equation}\label{eq:sumckeven}
 \sum_{k=4}^n c_k
 \leq 0.19\left( \left(\frac{n-2}{2}\right)! \right)^2+\frac{9}{4}\left(\left(\frac{n-2}{2}\right)! \right)^2
 =2.44\left(\left(\frac{n-2}{2}\right)! \right)^2.
 \end{equation}
 By Stirling's approximation, as $n\rightarrow\infty$,
 \begin{align*}
 2.44\left(\left(\frac{n-2}{2}\right)! \right)^2
 &=(1+o(1))2.44\left(\frac{n-2}{2e}\right)^{n-2}\cdot \pi(n-2)\\
 &=(1+o(1))2.44\pi \frac{n^{n-1}}{(2e)^n}(2e)^2\left(\frac{n-2}{n}\right)^{n-1}\\
 &=(1+o(1))2.44\pi \frac{n^{n-1}}{(2e)^n}4e^2\cdot \frac{1}{e^2}\\
 &=(1+o(1))9.76\pi \frac{n^{n-1}}{(2e)^n}\\
 &<10\pi \frac{n^{n-1}}{(2e)^n} &&\mbox{(for $n$ suff. large)}
 \end{align*}
  completing the proof of the lemma.\qed\bigskip

\begin{lemma}\label{le:1odd}
There exists $n_0\in\mathbb{Z}^{+}$ so that for every odd integer
$n\geq n_0$, if $G$ is a triangle-free graph on $n$ vertices, and
$x_1x_2\in E(G)$ with $\deg_G(x_2)\leq \frac{2}{5}n$,
then the number of cycles containing the edge $x_1x_2$
is at most $7.81\pi \frac{n^{n-1}}{(2e)^n}$.
\end{lemma}

\noindent{\bf Proof:}  The proof is similar to that of Lemma \ref{le:1even}.
 Let $G$ be a triangle-free graph on $n$ vertices,
and let $x_1x_2\in E(G)$, where $\deg(x_2)\leq \frac{2}{5}n$.
For each $k=4,\ldots,n$, let $c_k$ denote the number of cycles
of length $k$ that contain the edge $x_1x_2$.

For $3\leq i\leq \frac{n-5}{2}$,
an upper bound on $c_{2i-1}+c_{2i}$ is first calculated; to do so, count all possible
cycles of the form $x_1,x_2,\ldots,x_{2i-1}$ or $x_1,x_2,\ldots,x_{2i}$.
As in Lemma \ref{le:1even}, for each $j>1$, there are at most
$d_j=|N(x_j)\backslash\{x_1,\ldots,x_{j-1}\}|$ ways   to choose an $x_{j+1}$, and
  \begin{equation}\label{eq:odddjdj+1}
     d_jd_{j+1}\leq \left\lfloor \frac{n-j}{2}\right\rfloor \cdot \left\lceil \frac{n-j}{2}\right\rceil.
 \end{equation}
Using  (\ref{eq:odddjdj+1}) and the fact that $d_2\leq \frac{2}{5}n$,
the number of ways to choose vertices
 $x_3,x_4,\ldots,x_{2i-1}$ so that $x_1,x_2,x_3,x_4,\ldots, x_{2i-1}$ form a path
 is at most
 \begin{align}\label{eq:1starodd}
 \prod_{j=2}^{2i-2}d_j=d_2\prod_{j=3}^{2i-2}d_j
 \leq \frac{2}{5}n\prod_{j=1}^{i-2}(d_{2j+1}d_{2j+2})
& \leq \frac{2}{5}n\prod_{j=1}^{i-2}\left(\left\lfloor \frac{n-2j-1}{2}\right\rfloor \cdot \left\lceil \frac{n-2j-1}{2}\right\rceil\right)\notag\\
& =\frac{2}{5}n\prod_{j=1}^{i-2}\left(\frac{n-2j-1}{2}\right)^2.
 \end{align}
 If $x_{2i-1}x_1\in E(G)$, there is one cycle of length $2i-1$ and no cycles of length
 $2i$; if there is no such edge, there are no cycles of length $2i-1$ and
 at most $n-2i-1$ cycles of length $2i+1$. By these observations and (\ref{eq:1starodd}),
 \begin{equation}\label{eq:2starodd}
 c_{2i-1}+c_{2i}\leq (n-2i-1)\frac{2}{5}n\prod_{j=1}^{i-2}\left(\frac{n-2j-1}{2}\right)^2.
 \end{equation}
 To evaluate $\sum_{k=4}^nc_k$, separate the sum into three parts:
 \[ \sum_{k=4}^{n}c_k=c_4+ \sum_{k=5}^{n-5}c_k+ \sum_{k=n-4}^{n}c_k.\]
First,
\begin{equation}\label{eq:3starodd}
c_4\leq d_2d_3<n\cdot n=n^2.
\end{equation}
Next,
\begin{align}
 \sum_{k=5}^{n-5}c_k
 &=\sum_{i=3}^{(n-5)/2}(c_{2i-1}+c_{2i})\notag \\
 &\leq \sum_{i=3}^{(n-5)/2}\left[(n-2i-1) \frac{2}{5}n
          \prod_{j=1}^{i-2}\left(\frac{n-2j-1}{2}\right)^2\right]
                   &&\mbox{(by (\ref{eq:2starodd}))}\notag \\
  &= \frac{2}{5}n \sum_{i=3}^{(n-5)/2}\left[(n-2i-1)
          \prod_{j=1}^{i-2}\left(\frac{n-2j-1}{2}\right)^2\right]\notag\\
 &=\frac{2}{5}n\sum_{i=3}^{(n-5)/2}(n-2i-1)
          \left(\frac{\left(\frac{n-3}{2}\right)!}{\left(\frac{n-2i-1}{2}\right)!} \right)^2
 \notag \\
 &=\frac{2}{5}n\left( \left(\frac{n-3}{2}\right)! \right)^2
                             \sum_{j=3}^{\frac{n-5}{2}}\frac{2j}{(j!)^2}\notag\\
 &=\frac{2}{5}n\left( \left(\frac{n-3}{2}\right)! \right)^2
 \left(\sum_{j=1}^{\frac{n-5}{2}}\frac{2j}{(j!)^2} -\frac{2}{(1!)^2}-\frac{2\cdot 2}{(2!)^2}\right)\notag\\
 &< \frac{2}{5}n\left( \left(\frac{n-3}{2}\right)! \right)^2(3.19-3)    &&\mbox{(by (\ref{eq:I_12}))} \notag\\
 &=0.076n\left( \left(\frac{n-3}{2}\right)! \right)^2.\label{eq:4starodd}
 \end{align}

To count $\displaystyle\sum_{k=n-4}^nc_k$, note that  by (\ref{eq:3starodd}), there are at most
 \[\prod_{i=2}^{n-5}d_i=d_2\cdot\prod_{j=1}^{(n-7)/2}d_{2j+1}d_{2j+2}\leq \frac{2}{5}n\prod_{j=1}^{\frac{n-7}{2}}\left(\frac{n-2j-1}{2}\right)^2\]
 ways to choose a path $x_1,x_2,\ldots,x_{n-4}$, and by Lemma \ref{le:6vcxy},
 there are at most 9 ways to complete to a cycle (by paths that connect
 $x_{n-4}$ and $x_1$) in the graph $G\backslash\{x_1,\ldots,x_{n-5}\}$.
Therefore,
 \begin{equation}\label{eq:5starodd}
 \sum_{k=n-4}^nc_k\leq 9\cdot \frac{2}{5}n\prod_{j=1}^{\frac{n-7}{2}}\left(\frac{n-2j-1}{2}\right)^2
 =9\cdot \frac{2}{5}n\cdot \frac{\left(\left(\frac{n-3}{2}\right)! \right)^2}{(2!)^2}
 =\frac{9}{10}n\left(\left(\frac{n-3}{2}\right)! \right)^2.
 \end{equation}
Adding (\ref{eq:3starodd}), (\ref{eq:4starodd}), and (\ref{eq:5starodd}), as $n\rightarrow\infty$,
\begin{align}
\sum_{k=4}^nc_k
&\leq n^2 + 0.076n\left( \left(\frac{n-3}{2}\right)! \right)^2
                                     +\frac{9}{10}n\left( \left(\frac{n-3}{2}\right)! \right)^2 \notag\\
 &= n^2 + 0.976n\left( \left(\frac{n-3}{2}\right)! \right)^2 \label{eq:sumckodd}\\
 &= n^2 + (1+o(1))0.976n(n-3)\pi\left( \frac{n-3}{2e}\right)^{n-3}\notag\\
 &=     (1+o(1))0.976\pi n\cdot \frac{n^{n-2}}{(2e)^n}\left( \frac{n-3}{n}\right)^{n-2}(2e)^3\notag\\
 &=     (1+o(1))0.976\pi \cdot \frac{n^{n-1}}{(2e)^n}\frac{1}{e^3}(2e)^3\notag\\
 & = (1+o(1))7.808\pi \cdot \frac{n^{n-1}}{(2e)^n\notag}\\
&< 7.81 \pi\frac{n^{n-1}}{(2e)^n} &&\mbox{(for $n$ suff. large),}\notag
\end{align}
completing the proof.\qed\bigskip

\begin{lemma}\label{le:2ham}
Let $H$ be a triangle-free graph on $k$ vertices. Then $H$ has
at most $e^2\left(\frac{k}{2e}\right)^k$ hamiltonian cycles.
\end{lemma}
\noindent{\bf Proof:}  Let $x_1$ be the first vertex of a hamiltonian cycle.
For each $i\geq1$, there are at most $d_i=|N(x_i)\backslash\{x_1,\ldots,x_{i}\}|$
ways to choose a vertex $x_{i+1}$. Note that $N(x_i)\cap N(x_{i+1})=\emptyset$
because if the intersection contains some vertex $v$, then $v$, $x_i$, and $x_{i+1}$
form a triangle. Also,
\[|N(x_i)\backslash\{x_1,\ldots,x_{i}\}\cup N(x_{i+1})\backslash\{x_1,\ldots,x_{i+1} \}|
  \leq |V(H)\backslash\{x_1,\ldots,x_{i+1} \}|=k-i.\]
  Therefore,
  \begin{align*}
  d_i+d_{i+1}
        &=|N(x_i)\backslash\{x_1,\ldots,x_{i}\}|+|N(x_{i+1})\backslash\{x_1,\ldots,x_{i+1} \}|\\
        &=|N(x_i)\backslash\{x_1,\ldots,x_{i}\}\cup N(x_{i+1})\backslash\{x_1,\ldots,x_{i+1} \}|\\
        &\leq k-i,
\end{align*}
and thus $d_id_{i+1}\leq \left\lfloor\frac{k-i}{2}\right\rfloor \cdot\left\lceil\frac{k-i}{2}\right\rceil$.
\medskip

When $k$ is odd, the number of hamiltonian cycles is at most
\[\prod_{i=1}^{k-1}d_i=\prod_{j=1}^{\frac{k-1}{2}}d_{2j-1}d_{2j}
    \leq \prod_{j=1}^{\frac{k-1}{2}} \left\lfloor\frac{k-2j+1}{2}\right\rfloor
                         \cdot\left\lceil\frac{k-2j+1}{2}\right\rceil
\mspace{-2.0mu}    =\prod_{j=1}^{\frac{k-1}{2}} \left(\frac{k-2j+1}{2}\right)^2 \mspace{-10.0mu} =\left( \left(\frac{k-1}{2}\right)!\right)^2\]
and by  (\ref{eq:simplebdsn!}), this number is at most
\[\left( \frac{\left(\frac{k-1}{2}\right)^{\frac{k-1}{2}+\frac{1}{2}}}{e^{\frac{k-1}{2}-1}}\right)^2
\mspace{-3.0mu}=\frac{\left(\frac{k-1}{2} \right)^k}{e^{k-3}}\mspace{-1.0mu}=e^3\left(\frac{k-1}{k}\right)^k\left(\frac{k}{2e} \right)^k
\mspace{-4.0mu}= e^3\frac{1}{\left(1+\frac{1}{k-1}\right)^k}\left(\frac{k}{2e} \right)^k
\mspace{-2.0mu}\leq e^2\left(\frac{k}{2e}\right)^k,\]
completing the proof for odd $k$.\medskip

When $k$ is even, similarly obtain
\begin{align*}
\prod_{i=1}^{k-1}d_i&=\left(\prod_{j=1}^{\frac{k-2}{2}}d_{2j-1}d_{2j}\right)\cdot d_{k-1}
    \leq \left(\prod_{j=1}^{\frac{k-2}{2}} \left\lfloor\frac{k-2j+1}{2}\right\rfloor
                         \cdot\left\lceil\frac{k-2j+1}{2}\right\rceil\right)\cdot1\\
   & =\prod_{j=1}^{\frac{k-1}{2}} \left(\frac{k-2j}{2}\right) \left(\frac{k-2j+2}{2}\right)
 \mspace{-1.0mu}   =\frac{k}{2}\left( \left(\frac{k-2}{2}\right)!\right)^2
 \mspace{-4.0mu}   \leq\frac{k}{2}\left( \frac{\left(\frac{k-2}{2}\right)^{\frac{k-2}{2}
                +\frac{1}{2}}}{e^{\frac{k-2}{2}}}\right)^2\\
&=k\frac{(k-2)^{k-1}}{e^{k-4}2^k}
=e^4\left(\frac{k-2}{k}\right)^{k-1}\left(\frac{k}{2e} \right)^k
\mspace{-4.0mu}= e^4\frac{1}{\left(1+\frac{2}{k-2}\right)^{k-1}} \left(\frac{k}{2e} \right)^k
\mspace{-4.0mu}\leq e^2\left(\frac{k}{2e}\right)^k,
\end{align*}
completing the proof for even $k$, and hence for the lemma.\qed
\bigskip

\section{Main theorems}
In Theorem \ref{th:main},  Conjecture \ref{con:main} is proved for sufficiently large $n$.
Then in Theorem \ref{th:mainngeq141}, a lower bound on such $n$ is given.
\begin{theorem}\label{th:main}
There exists $n_0\in\mathbb{Z}^{+}$ so that for any $n\geq n_0$, the triangle-free
graph on $n$ vertices with the largest number of cycles is $K_{\lfloor n/2\rfloor, \lceil n/2\rceil}$.
\end{theorem}
\noindent{\bf Proof:} Let $G$ be a triangle-free graph on $n$ vertices. It is first shown
that if $G$ contains a vertex of small degree, then $G$ has far fewer cycles than does
$K_{\lfloor n/2\rfloor, \lceil n/2\rceil}$.

Let $x\in V(G)$, and assume that $\deg(x)\leq \frac{2}{5}n$.
Cycles in $G$ are counted according to whether or not they contain $x$.

The number of cycles not containing $x$:
Any cycle in $G-x$ is a hamiltonian cycle for some subgraph,
and so the number of cycles in $G$ not containing $x$ is loosely bounded above
by
\begin{align}
\sum_{Y\subseteq V(G) \backslash x} &(\mbox{number of ham. cycles in $G[Y]$})\\
&\leq \sum_{k=4}^{n-1}\binom{n-1}{k}e^2\left( \frac{k}{2e}\right)^k
                                                     && (\mbox{by Lemma \ref{le:2ham}})\notag\\
& <  e^2\sum_{k=4}^{n-1}\binom{n-1}{k}\left( \frac{n-1}{2e}\right)^k\notag\\
& < e^2\left(1+\frac{n-1}{2e}\right)^{n-1}\notag\\
& = e^2\left(\frac{n+2e-1}{2e}\right)^{n-1}\notag\\
& = e^2\left(\frac{n}{2e}\right)^{n-1}\left(\frac{n+2e-1}{n}\right)^{n-1}\notag\\
& < e^2\left(\frac{n}{2e}\right)^{n-1}\left(1+\frac{2e-1}{n}\right)^{n}\notag\\
& \leq  e^2\left(\frac{n}{2e}\right)^{n-1}e^{2e-1}\notag\\
& =\frac{2e^{2e+2}}{n}\left(\frac{n}{2e}\right)^n. \label{eq:cyclesnox}
\end{align}

The number of cycles containing $x$:
Each cycle $C$ containing $x$ has exactly two edges (in $C$) incident with $x$,
and so the number of cycles containing $x$ is
\begin{equation}\label{eq:cincwxy}
\frac{1}{2}\sum_{y\in N(x)}(\mbox{number of cycles containing $xy$}).
\end{equation}

By Lemma \ref{le:1even}, for even $n$, the expression (\ref{eq:cincwxy}) is at most
\[\frac{1}{2}\cdot \frac{2}{5}n\cdot 10\pi\frac{n^{n-1}}{(2e)^n}=2\pi \left(\frac{n}{2e}\right)^n.\]
In this case, for $n$ sufficiently large, the total number of cycles in $G$ is at most
\[2\pi\left(\frac{n}{2e}\right)^n+\frac{2e^{2e+2}}{n}\left(\frac{n}{2e}\right)^n
=\left(2\pi +\frac{2e^{2e+2}}{n}\right)\left(\frac{n}{2e}\right)^n\leq 2.01\pi \left(\frac{n}{2e}\right)^n.
\]
However, by (\ref{eq:cTn2lowerbd}), the number of cycles in $K_{\lfloor n/2\rfloor, \lceil n/2\rceil}$ is (for $n$ even) at least  $2.27958\pi\left(\frac{n}{2e}\right)^n$.\medskip

Let $n$ be odd; then by Lemma \ref{le:1odd}, the expression (\ref{eq:cincwxy}) is
at most
\begin{equation}\label{eq:oddcxy}
\frac{1}{2}\cdot \frac{2}{5}n\cdot 7.81\pi\frac{n^{n-1}}{(2e)^n}
=1.562\pi\left(\frac{n}{2e}\right)^n.
\end{equation}
Thus, for odd $n$ sufficiently large,
by (\ref{eq:oddcxy}) and (\ref{eq:cyclesnox}) the total number of cycles in $G$ is at most
\[1.562\pi\left(\frac{n}{2e}\right)^n+\frac{2e^{2e+2}}{n}\left(\frac{n}{2e}\right)^n
                  \leq 1.57\pi\left(\frac{n}{2e}\right)^n.\]
By (\ref{eq:cTn2lowerbd}) in Theorem \ref{th:cTn2approx}, the number of cycles in
$K_{\lfloor n/2\rfloor, \lceil n/2\rceil}$ for $n$ odd is at least  $1.5906\pi\left(\frac{n}{2e}\right)^n$.

In both the even and odd case, if $G$ contains a vertex of degree at most $\frac{2}{5}n$,
 then $G$ has far fewer cycles than does $K_{\lfloor n/2\rfloor, \lceil n/2\rceil}$.

So assume that $\delta(G)>\frac{2}{5}n$. Then by Theorem \ref{th:andr64}, $G$ is
bipartite. By Lemma \ref{le:maxcbipartite}, the number of cycles in $G$ is maximized
by $K_{\lfloor n/2\rfloor, \lceil n/2\rceil}$.
\qed\bigskip

\begin{theorem} \label{th:mainngeq141}
The statement of Theorem \ref{th:main} with $n_0=141$ is true.
\end{theorem}
\noindent{\bf Proof:} To show that $n_0$ works, further estimations on
$c(K_{\lfloor n/2\rfloor, \lceil n/2\rceil})$ are needed for $n\geq 141$.
Both when $n$ is even and when $n$ is odd,  (\ref{eq:cT2nfracdotan}) holds (but
the expression for $a_\ell$ changes). Since each (one for odd, one for even) sequence of $a_\ell$s are non-increasing
 for $n\geq 140$,
 \begin{align}
 c(K_{\lfloor n/2\rfloor, \lceil n/2\rceil})
 & \leq \frac{\lfloor n/2\rfloor ! \lceil n/2\rceil !}{2\lfloor n/2 \rfloor} \cdot
     \begin{cases}
         a_{71} &\mbox{for $n$ even}\\
         a_{70} &\mbox{for $n$ odd}
  \end{cases}  \notag \\
  \label{eq:cTn2upper141}
 & \leq \frac{\lfloor n/2\rfloor ! \lceil n/2\rceil !}{2\lfloor n/2 \rfloor}
 \begin{cases}
2.302786 &\mbox{for $n$ even}\\
1.60067 &\mbox{for $n$ odd}
  \end{cases}.
\end{align}
(The values of $a_{70}$ and $a_{71}$ were calculated by computer.)
With these estimates in hand, now Theorem \ref{th:main} is proved with $n_0=141$.
Let $G$ be a triangle-free graph on $n\geq 141$ vertices.
Without loss of generality, assume that there is a vertex
of degree at most $\frac{2}{5}n$ (since otherwise, the theorem is proved by
Theorem \ref{th:andr64} and Lemma \ref{le:maxcbipartite}). In the following
calculations, bounds given in (\ref{eq:simplebdsn!}) and Theorem \ref{th:cTn2approx}
are used freely.

Case 1: Let $n\geq 141$ be odd.
By (\ref{eq:sumckodd}) from
the proof of Lemma \ref{le:1odd}, the number of
cycles passing through an edge $xy$ in $G$ is at most
$ n^2 + 0.976n\left( \left(\frac{n-3}{2}\right)! \right)^2$.
Then the number of cycles in $G$ is bounded by
\begin{align*}
c(G)& \leq\frac{1}{2}\cdot \frac{2}{5}n
        \cdot   \mspace{-2.0mu}\left[n^2 + 0.976n\left( \left(\frac{n-3}{2}\right)! \right)^2\right]
        + \frac{2e^{2e+2}}{n}\left(\frac{n}{2e}\right)^n\\
   &=\frac{\frac{n-1}{2}!\frac{n+1}{2}!}{n-1}\cdot I_1(2) \cdot \left(\frac{\frac{n}{5} \left[n^2 + 0.976n\left( \left(\frac{n-3}{2}\right)! \right)^2\right](n-1)}{\frac{n-1}{2}!\frac{n+1}{2}!\cdot I_1(2)} \right) \mspace{-2.0mu}
   +I_1(2)\cdot \pi\left(\frac{n}{2e}\right)^n\left(  \frac{2e^{2e+2}}{n\pi I_1(2)}  \right)\\
& \leq   c(K_{\lfloor n/2\rfloor, \lceil n/2\rceil})\cdot
 \left(10^{-10}+\frac{8}{5}\cdot (0.976)\left(\frac{n^2}{n^2-1}\right) +\frac{2e^{2e+2}}{n}\right)
 \cdot \frac{1}{I_1(2)}\\
& \leq   c(K_{\lfloor n/2\rfloor, \lceil n/2\rceil})\cdot 6.
\end{align*}

Case 2: Let $n$ be even and $n\geq 142$. Then by (\ref{eq:sumckeven}), the proof
of Theorem \ref{th:main}, and by the result in Case 1,
\begin{align*}
c(G) &\leq \frac{1}{2}\cdot \frac{2}{5}n\cdot 2.44\left(\left(\frac{n-2}{2}\right)! \right)^2+
  6\cdot c(K_{\lfloor n/2\rfloor, \lceil n/2\rceil})\\
&  = \frac{\frac{n}{5}2.44\left(\left(\frac{n-2}{2}\right)! \right)^2}{\frac{\lfloor n/2\rfloor ! \lceil n/2 \rceil !}{2\lfloor n/2\rfloor }\cdot I_0(2)}
        \cdot \frac{\lfloor n/2\rfloor ! \lceil n/2 \rceil !}{2\lfloor n/2\rfloor}\cdot I_0(2)
     + c(K_{\lfloor n/2\rfloor, \lceil n/2\rceil})\frac{6\cdot c(K_{\lfloor \frac{n-1}{2}\rfloor, \lceil \frac{n-1}{2}\rceil})}{c(K_{\lfloor n/2\rfloor, \lceil n/2\rceil})}\\
     &\leq c(K_{\lfloor n/2\rfloor, \lceil n/2\rceil})\cdot \left( \frac{\frac{4}{5}\cdot 2.44}{I_0(2)}+\frac{6\cdot 1.60067 \cdot  \frac{\lfloor \frac{n-1}{2}\rfloor ! \lceil \frac{n-1}{2} \rceil !}{2\lfloor \frac{n-1}{2}\rfloor}}{I_0(2)  \frac{\lfloor n/2\rfloor ! \lceil n/2 \rceil !}{2\lfloor n/2\rfloor}} \right)\\
 &=  c(K_{\lfloor n/2\rfloor, \lceil n/2\rceil})\left(\frac{\frac{4}{5}\cdot 2.44}{I_0(2)}+
 \frac{6\cdot 1.60067}{I_0(2)}\cdot \frac{2}{n}  \right)\\
 & \leq c(K_{\lfloor n/2\rfloor, \lceil n/2\rceil})  \mbox{\ \ \ (for $n\geq 142$)}.
\end{align*}
Returning to the case when $n$ is odd,
\begin{align*}
c(G) & \leq c(K_{\lfloor n/2\rfloor, \lceil n/2\rceil})\cdot
   \left(\frac{10^{-10}+\frac{8}{5}\cdot (0.976)\left(\frac{n^2}{n^2-1}\right)}{I_1(2)}
   +\frac{2.302786 \cdot
   \frac{\lfloor \frac{n-1}{2}\rfloor ! \lceil \frac{n-1}{2} \rceil !}{2\lfloor \frac{n-1}{2} \rfloor}}{I_1(2)
      \cdot \frac{\lfloor n/2\rfloor ! \lceil n/2 \rceil !}{2\lfloor n/2\rfloor}}\right)\\
  &\leq c(K_{\lfloor n/2\rfloor, \lceil n/2\rceil})\cdot
   \left(\frac{10^{-10}+\frac{8}{5}\cdot (0.976)\left(\frac{n^2}{n^2-1}\right)}{I_1(2)}
   +\frac{2.302786 }{I_1(2)  \cdot (n+1)}\right)\\
 &\leq c(K_{\lfloor n/2\rfloor, \lceil n/2\rceil})\cdot 0.9947\\
 & <c(K_{\lfloor n/2\rfloor, \lceil n/2\rceil}).
\end{align*}
This completes the proof of the theorem for $n\geq 141$. \qed\bigskip

\section{Concluding remarks}
%In the proof of Theorem \ref{th:main}, it is shown that for large $n$,
%if a triangle-free graph $G$ on $n$ vertices is not
%bipartite, then $G$ has significantly fewer cycles than does
%$K_{\left \lceil \frac{n}{2}\right \rceil, \left \lfloor \frac{n}{2}\right \rfloor}$.

A few natural extensions of Theorem \ref{th:main} and Conjecture \ref{con:main} can be considered.  For example, for $k>1$, what $C_{2k+1}$-free graphs have the most number
of cycles?  It is well known (see, {\it e.g}, \cite[p.~150]{BB:79}, \cite{FuGu:14}, or
\cite{Simo:72}) that for $n$ large enough, the unique $C_{2k+1}$-free $n$-vertex
graph with the maximum
number of edges is $K_{\left \lceil \frac{n}{2}\right \rceil ,\left \lfloor \frac{n}{2}\right \rfloor}$.
Also, Balister, Bollob\'as, Riordan, and Schelp \cite{BBRS:03} showed that for a fixed $n$ and $k$ and $\frac{n}{2}<\Delta<n-k$, any maximal $C_{2k+1}$-free
graph with maximum degree $\Delta$ on $n$ vertices is the complete bipartite graph
$K_{\Delta, n-\Delta}$.  Since any $C_{2k+1}$-free graph with a maximum number of
cycles is also edge-maximal, it might then seem reasonable to
pose the following:
\begin{conjecture}\label{con:oddcyclefree}
For any $k>1$, if an $n$-vertex graph $C_{2k+1}$-free graph has the maximum number
of cycles, then $G=K_{\left \lceil \frac{n}{2}\right \rceil ,\left \lfloor \frac{n}{2}\right \rfloor}$.
\end{conjecture}
In support of Conjecture \ref{con:oddcyclefree}, by duplicating the proofs in this paper, it can be shown that the order of magnitude for edges in a $C_{2k+1}$-free graph is correct:
\begin{theorem} \label{th:oddcyclefreeorder}
For any $k\geq 2$, there exists a constant $\alpha_{2k+1}\geq 1$ so that
 if $G$ is a $C_{2k+1}$-free graph on $n$ vertices with the maximum number of cycles,
 then $c(G)\leq \alpha_{2k+1}\cdot c(K_{\left \lceil \frac{n}{2}\right \rceil ,\left \lfloor \frac{n}{2}\right \rfloor})$.
 \end{theorem}
 With considerably more work, it seems feasable that by refining methods in this paper,
when $k=2$, the constant $\alpha_5$ can be reduced to 1. However, for $k>2$, our proof of
 Theorem \ref{th:oddcyclefreeorder} yields a  constant $\alpha_{2k+1}$ that
  grows exponentially in $k$, so it seems that new techniques are required to
 settle Conjecture \ref{con:oddcyclefree} for general $k$.  In any case, for $k>2$, to prove the
 general structure of a cycle-maximal $C_{2k+1}$-free graph seems beyond our reach
 at this moment.

One might entertain other questions related to Conjecture \ref{con:main}. For example:
\begin{question}\label{qu:cyclesgirth}
What is the maximum number of cycles in a graph on $n$ vertices with girth at least $g$?
\end{question}
The case $g=3$ is trivial and this paper addresses this question for $g=4$;
 however, there seems to be little known for $g\geq5$.  Another question
that might be interesting is:
\begin{question} \label{qu:cyclescliquefree}
For $k\geq 4$  what is the  maximum number of cycles in a $K_k$-free graph on $n$ vertices?
Could it be that the cycle-maximal $K_k$-free graphs are indeed Tur\'an graphs?
\end{question}

A type of stability result also follows from the techniques given in this paper.
Theorem \ref{th:main} shows that among all triangle-free graphs with
$n$ vertices and $m=\left \lfloor \frac{n^2}{4}\right \rfloor$ edges,
$K_{\left \lceil \frac{n}{2}\right \rceil ,\left \lfloor \frac{n}{2}\right \rfloor}$ has
the most number of cycles.
Let $\ell=o(n)$, and set $m=\left \lfloor \frac{n^2}{4}\right \rfloor-\ell$.
If $G$ has $n$ vertices and $m$ edges, and has the most number of cycles
among all triangle-free $n$-vertex graphs with $m$ edges, then same argument as in the proof of Theorem \ref{th:main} implies that $G$ is bipartite.
By the maximality of the number of cycles, one can show that $G$ is a
subgraph of $K_{\left \lceil \frac{n}{2}\right \rceil ,\left \lfloor \frac{n}{2}\right \rfloor}$.

% Denote by $MC(n,m)$ the maximal number of cycles in a triangle-free graph on $n$ vertices that has at most $m$ edges. Then Theorem \ref{th:mainngeq141}  states that for $n\geq 141$,
%\[MC(n, \left \lfloor{ \frac{n^2}{4}}\right \rfloor )
%           =c(K_{\left \lfloor{ \frac{n}{2}}\right \rfloor ,\left \lceil{\frac{n}{2}}\right \rceil} ).\]
%With the same approach as in Theorem \ref{th:main}, we can prove that
%         \[MC(n,\left \lfloor{ \frac{n^2}{4}}\right \rfloor -o(n))=c(G),\]
%where $G$ is some subgraph of  $K_{\left \lfloor{ \frac{n}{2}}\right \rfloor ,\left \lceil{\frac{n}{2}}\right \rceil}$.
%This can be considered as a stability version of Theorem \ref{th:main}.
%
For $14\leq n \leq 140$, Conjecture \ref{con:main} remains open.  With a bit more care, it appears that  with the techniques in this paper, one might be able to
prove Conjecture \ref{con:main} for the even $n$ to $n\geq 100$ or so, but the techniques
used here do not seem to leave much room for  the odd $n$.
Skala \cite{Skala:15} has suggested that
Lemma \ref{le:6vcxy} might be proved for graphs with slightly more vertices,
and such an improvement might yield modest improvements for the bound on $n$
for which Theorem \ref{th:main} holds.

Finally, during preparation of this paper, Alex Scott (Oxford) has informed us that
he and his students have been working on Conjecture \ref{con:main}; he reports
that they had some success using the Regularity Lemma (which would give a
value for $n_0$ much larger than 141) , and that they are
working on a different proof of a more general result.

\end{document}